\documentclass{ifacconf}
\usepackage{mathtools}
\usepackage{amstext}
\usepackage{amssymb}

\usepackage{bm}
\usepackage{color}

\renewcommand{\S}{\mathcal{S}}
\newcommand{\I}{\mathcal{I}}
\newcommand{\G}{\mathcal{G}}
\newcommand{\V}{\mathcal{V}}
\newcommand{\E}{\mathcal{E}}
\newcommand{\x}{\mathbf{x}}
\renewcommand{\u}{\mathbf{u}}
\newcommand{\w}{\mathbf{w}}
\renewcommand{\c}{\mathbf{c}}
\newcommand{\boldb}{\bm{\beta}}
\newcommand{\boldg}{\bm{\gamma}}
\newcommand{\tL}{\mathcal{L}}
\usepackage{flushend}

\usepackage{graphicx}      
\usepackage[labelformat=simple]{subcaption}

\usepackage{natbib}
\begin{document}
\begin{frontmatter}

\title{Sparse estimation of Laplacian eigenvalues in multiagent networks} 

\thanks[footnoteinfo]{This work was supported, in part, by the National Science Foundation, grant CAREER-ECCS-1651433.}

\author[First]{Mikhail Hayhoe}
\author[First]{Francisco Barreras}
\author[First]{Victor M. Preciado} 

\address[First]{University of Pennsylvania, 
   Philadelphia, PA 19143 USA (e-mail: mhayhoe@seas.upenn.edu, fbarrer@sas.upenn.edu, preciado@seas.upenn.edu).}

\begin{abstract}                
We propose a method to efficiently estimate the Laplacian eigenvalues of an arbitrary, unknown network of interacting dynamical agents. The inputs to our estimation algorithm are measurements about the evolution of a collection of agents (potentially one) during a finite time horizon; notably, we do not require knowledge of which agents are contributing to our measurements. We propose a scalable algorithm to exactly recover a subset of the Laplacian eigenvalues from these measurements. These eigenvalues correspond directly to those Laplacian modes that are observable from our measurements. We show how our technique can be applied to networks of multiagent systems with arbitrary dynamics in both continuous- and discrete-time. Finally, we illustrate our results with numerical simulations.
\end{abstract}


\end{frontmatter}

\section{Introduction}

The spectrum of the Laplacian matrix describing a network of interacting dynamical agents provides a wealth of global information about the network structure and function; see, e.g., \cite{Fiedler73, mohar91, merris94, Chung97, mesbahi10, Bullo19}, and references therein.
For example, the Laplacian spectrum finds applications in multiagent coordination problems as in \cite{jadbabaie03, olfati07}, synchronization of oscillators in \cite{pecora98, dorfler13}, neuroscience such as \cite{Becker18}, biology as in \cite{Palsson06}, as well as several graph-theoretical problems, such as finding cuts (see \cite{shi00}) or communities (see \cite{Luxburg07}) in graphs, among many others, as illustrated in \cite{mohar97}.

Due to its practical importance, numerous methods have been proposed to estimate the Laplacian eigenvalues of a network of dynamical agents.
For example, \cite{kempe08} proposed a distributed algorithm based on orthogonal iteration (see \cite{golub13}) for computing higher-dimensional invariant subspaces.
In the control literature, \cite{Franceschelli13} define local interaction rules between agents such that the network response is a superposition of sinusoids oscillating at frequencies related to the Laplacian eigenvalues; however, this approach imposes a particular dynamics on the agents in the network, which is unrealistic in many scenarios.
\cite{Aragues14} proposed a distributed algorithm based on the power iteration for computing upper and lower bounds on the algebraic connectivity (i.e., the second smallest Laplacian eigenvalue).
An approach by \cite{Kibangou15} uses consensus optimization to deduce the spectrum of the Laplacian, but this requires a consensus algorithm to be run on the network separately from the dynamics.
Using the Koopman operator, it has been shown that the spectrum of the Laplacian may be recovered using sparse local measurements, see \cite{Mauroy17,Mesbahi19}; 
unfortunately, these methods require the system to be reset to known initial conditions multiple times.
\cite{Leonardos20} proposed a distributed continuous-time dynamics over manifolds to compute the largest (or smallest) eigenvalues and eigenvectors of a graph.

We find in the literature several works more closely related to the techniques used in this paper. For example, we find a classical result in linear algebra, referred to as the \emph{Newton-Girard equations} (see, e.g., \cite{Herstein06}) which allows us to recover eigenvalues by analyzing symmetric polynomials of the traces of powers of the matrix. In a similar line of work, \cite{Preciado13-2} used the spectral moments of a graph, computed from local structural information, to compute bounds on spectral properties of practical importance. Another related method uses tools from probability theory to approximate the spectrum of a graph by counting the number of walks of length $k$ and then solving the classical moment problem, as in \cite{Preciado13, Barreras19}. The latter approach requires only local measurements of walks, but provides only bounds on the support of the eigenvalue spectrum. Apart from estimating the eigenvalues of a graph, we also find works aiming to reconstruct the whole graph structure from the dynamics, such as \cite{Shahrampour13, Shahrampour14}.

In this paper, we propose an approach which uses only a single sequence of measurements from the dynamics of a multiagent system, without prior knowledge of the network topology or initial condition. These measurements can be taken locally from a single agent, or from any linear combination of agents' outputs; in any case, we do not require knowledge of which agents are contributing to our measurements. From this single sequence of measurements we provide an efficient algorithm to estimate a subset of Laplacian eigenvalues. In particular, we estimate the eigenvalues corresponding to the observable modes of the Laplacian dynamics. In comparison to other approaches, our algorithm requires no parameter tuning, as the performance does not depend on any parameters. These techniques are applied in both discrete- and continuous-time to networks of single integrators, as well as more general multi-agent systems.

The remainder of this paper is structured as follows. We outline background and notation for our problem in Section~\ref{sec:bg}. In Section~\ref{sec:dt} we present our results for discrete-time systems, and in Section~\ref{sec:ct} we describe our results in the continuous-time case. Section~\ref{sec:sims} illustrates our results via simulations in a variety of systems, and Section~\ref{sec:conc} concludes the paper.

\section{Background and Notation}\label{sec:bg}

\begin{table}[ht]
    \centering
    \begin{tabular}{c|l}
        Symbol & Meaning \\
        \hline
        $I_n$ & $n\times n$ identity matrix \\
        $\mathbb{R}$ & set of real numbers \\
        $\mathbf{e}_i$ & $i$-th vector in the canonical basis of $\mathbb{R}^n$ \\
        $\V$ & node set, $\V = \{1,\ldots,n\}$ \\
        $\E$ & edge set, $\E \subseteq \V\times\V$ \\
        $\G = (\V,\E) $ & graph with node set $\V$ and edge set $\E$\\
        $\otimes$ & Kronecker product \\
        $\sigma(X) \coloneqq \{\lambda_i\}_{i=1}^n$ & eigenvalue spectrum of matrix $X$ \\
        $\chi(X)$ & characteristic polynomial of matrix $X$ \\
        $\delta(\cdot)$ & Dirac delta function \\
        $\mu_X \coloneqq \sum_{i=1}^n\delta(\lambda_i)$ & spectral measure of matrix $X$ \\
        $G = G(\G)$ & adjacency matrix of $\G$, $[G]_{ij} = 1 \Leftrightarrow \{i,j\}\in\E$ \\
        $D = D(\G)$ & degree matrix of $\G$, $[D]_{ii} = \sum_{j=1}^n [G]_{ij}$ \\
        $L = L(\G)$ & Laplacian matrix of $\G$, $L = D - G$ \\
        $\tL = \tL(\G)$ & normalized Laplacian of $\G$, $\tL =D^{-1}G$
    \end{tabular}
    \label{tab:notation}
\end{table}

Throughout this paper we use lower-case letters for scalars, lower-case bold letters for vectors, upper-case letters for matrices, and calligraphic letters for sets.

An undirected graph $\G = (\V,\E)$ has node set $\V$ and edge set $\E$, where $\{i,j\}\in\E$ means nodes $i$ and $j$ are connected. In this paper, we assume $\G$ is a simple graph. 


\section{Spectral Estimation for Discrete-Time Dynamics}\label{sec:dt}

We begin our exploration by considering the discrete-time (DT) dynamics of a network of single integrators. In this context, we will present a methodology to estimate the eigenvalues of the Laplacian matrix from a finite sequence of  measurements of our system. In Subsection~\ref{subsec:dt_multi}, we will extend this result to more general discrete-time agent dynamics, and will consider the continuous-time (CT) case in Section~\ref{sec:ct}.

\subsection{Network of Discrete-Time Single Integrators}\label{subsec:dt_single}

Consider the discrete-time dynamics of a collection of single integrators,
\begin{align}\begin{split}\label{eq:dt_integrator}
\mathbf{x}\left[k+1\right] & = \tL\mathbf{x}\left[k\right],\;\mathbf{x}\left[0\right]=\mathbf{x}_{0},\\
y\left[k\right] & =\mathbf{c}^{\intercal}\mathbf{x}\left[k\right],
\end{split}
\end{align}
where $\tL \coloneqq D^{-1}G$ is the normalized Laplacian matrix, $k\in\mathcal{\mathbb{N}}$, and $\mathbf{c},\mathbf{x}_{0}$
are arbitrary (possibly unknown) vectors in $\mathbb{R}^{n}$. For example, we may have
$\mathbf{c}=\mathbf{e}_{i}$ when we only observe the state of agent
$i$, or $\mathbf{c}=\sum_{i\in\mathcal{S}\subseteq\mathcal{V}}\beta_i\mathbf{e}_{i}$
when we observe the weighted sum of the states of a subset $\mathcal{S}$ of
agents. Thus, the evolution of the output measurement is
\[
y\left[k\right]=\mathbf{c}^{\intercal}\tL^{k}\mathbf{x}_{0}.
\]
In what follows we propose an algorithm to recover the eigenvalues of the normalized Laplacian $\tL$ from the output sequence $y[0],y[1],\ldots,y[2n-1]$.

Let $\mathbf{u}_{i}$ and $\mathbf{w}_{i}$ be the $i$-th right and left eigenvectors of $\tL$, respectively. Since $\tL$ is always diagonalizable with real eigenvalues (see \cite{Chung97}), we have that $\tL=U\Lambda W$, where $\Lambda \coloneqq \text{diag}(\lambda_1,\ldots,\lambda_n)$, $U \coloneqq [\u_1,\ldots,\u_n]$, and  $W \coloneqq [\w_1^\intercal;\cdots;\w_n^\intercal] = U^{-1}$; hence,
\begin{align*}
y\left[k\right] & =\left(\mathbf{c}^{\intercal}U\right)\Lambda^{k}\left(W\mathbf{x}_{0}\right) =\sum_{i=1}^{n}\omega_{i}\lambda_{i}^{k},
\end{align*}
where the $\lambda_{i}$ are real, and the weights $\omega_i$ are given by
\begin{equation}
\omega_{i} \coloneqq \left[\mathbf{c}^{\intercal}U\right]_{i}\left[W\mathbf{x}_{0}\right]_{i}=\mathbf{c}^{\intercal}\mathbf{u}_{i}\mathbf{w}_{i}^{\intercal}\mathbf{x}_{0}.\label{eq:weights}
\end{equation}

Define the following \emph{signed} Borel measure on $\mathbb{R}$:
\[
\mu_{\tL}(z) \coloneqq \sum_{i=1}^{n}\omega_{i}\delta\left(z-\lambda_{i}\right),
\]
which we refer to as the \emph{spectral measure of $\tL$}.
From (\ref{eq:weights}), we see that it is possible for $\omega_{i}=0$ whenever $\mathbf{c}^{\intercal}\mathbf{u}_{i}=0$ or $\mathbf{w}_{i}^{\intercal}\mathbf{x}_{0}=0$. Notice that if $\x_0$ is randomly generated, then almost surely $\w_i^\intercal\x_0 \neq 0$; hence, $\omega_i = 0$ for those eigenvalues $\lambda_i$ for which $\mathbf{c}^{\intercal}\mathbf{u}_{i}=0$. Therefore, those eigenvalues corresponding to unobservable eigenmodes of the Laplacian dynamics, according to the Popov-Belevitch-Hautus (PBH) test (see \cite{hespanha18}), will have $\omega_i=0$ and it will be impossible to recover them from our observations. Moreover, for some (possibly deterministic) initial condition $\x_0$, there are other (observable) eigenvalues that our method will not be able to recover. In particular, it may be that for some repeated eigenvalue $\lambda_{i}$, we have $\sum_{j:\lambda_{j}=\lambda_{i}}\omega_{j}=0$. Hence, the support of $\mu_{\tL}$ is the set
\[
\S_{\mu_{\tL}}\coloneqq\left\{ \lambda_{i}\in\sigma\left(\tL\right)\colon\sum_{j:\lambda_{j}=\lambda_{i}}\omega_{j}\neq0\right\}.
\]
Notice that, for a random initial condition $\x_0$, the set $\S_{\mu_{\tL}}$ almost surely coincides with the set of eigenvalues corresponding to observable eigenmodes in the PBH test.
In what follows, we state that the eigenvalues in this support set are those that can be recovered by any algorithm using these measurements, as demonstrated below.

\begin{lem}\label{lem:conds}
The eigenvalues which may be recovered from the sequence of measurements $(y[k])_{k=0}^s$, for any finite $s$, are exactly those in $\S_{\mu_{\tL}}$.
\end{lem}

\begin{pf}
Notice that
\begin{align*}
    y[k] &= \sum_{i \in \S_{\mu_{\tL}}} \omega_i\lambda_i^k + \sum_{i \not\in \S_{\mu_{\tL}}} \omega_i\lambda_i^k \\
    &= \sum_{i \in \S_{\mu_{\tL}}} \omega_i\lambda_i^k + \sum_{i \not\in \S_{\mu_{\tL}}: i \leq j~\forall j \text{ s.t. } \lambda_j=\lambda_i} \lambda_i^k \sum_{j:\lambda_j=\lambda_i}\omega_j.
\end{align*}
By definition of $\S_{\mu_{\tL}}$, the term $\sum_{j:\lambda_{j}=\lambda_{i}}\omega_{j}$ is zero. Hence, any eigenvalue $\lambda_i \not\in \S_{\mu_{\tL}}$ will never appear in any observation $y[k]$; therefore, this eigenvalue may not be recovered from any finite sequence of measurements.\qed
\end{pf}

In our algorithm, we use some tools from probability theory, introduced below. The $k$-th moment $m_k$ of the measure $\mu_{\tL}$ is defined by
\begin{align}\label{eq:moments}
m_{k} &\coloneqq \int_{\mathbb{R}}z^{k}d\mu_{\tL}\left(z\right) \cr
&= \int_{\mathbb{R}}z^{k}\sum_{i=1}^{n}\omega_{i}\delta\left(z-\lambda_{i}\right)dz \cr         
    &= \sum_{i=1}^{n}\omega_{i}\lambda_{i}^{k}.
\end{align}
Therefore $y\left[k\right]=m_{k}$, i.e., the $k$-th observation from our system is also the $k$-th moment of the spectral measure $\mu_{\tL}$. In what follows, we will propose a computationally efficient methodology to recover the support of $\mu_{\tL}$ using the sequence $\left(y\left[k\right]\right)_{k=0}^{2n-1} = \left(m_k\right)_{k=0}^{2n-1}$. Towards that goal, we define the Hankel matrix of moments
\begin{equation}
H \coloneqq \left[\begin{array}{cccc}
m_{0} & m_{1} & \cdots & m_{n-1}\\
m_{1} & m_{2} & \cdots & m_{n}\\
\vdots & \vdots & \ddots & \vdots\\
m_{n-1} & m_{n} & \cdots & m_{2n-2}
\end{array}\right].\label{eq:hankel}
\end{equation}
The following result relates the rank of this Hankel matrix to the cardinality of the support of $\mu_{\tL}$.

\begin{lem}
\label{lem:rank}The rank of $H$ in \eqref{eq:hankel} satisfies
\[
\text{rk}(H) = |\S_{\mu_{\tL}}|.
\]
\end{lem}

\begin{pf}
Let us define $\I \coloneqq \left\{ i\in\{1,\ldots,n\} : \lambda_i \in S_{\mu_{\tL}}\right\}$ and $\bar{\omega}_{i} \coloneqq \sum_{j:\lambda_{j}=\lambda_{i}}\omega_{j}$ for $i \in \I$.
Thus, we have $m_{k}=\sum_{i\in\I}\bar{\omega}_{i}\lambda_{i}$.
Now let $\mathbf{v}_{i}\coloneqq\left[1,\lambda_{i},\lambda_{i}^{2}\ldots,\lambda_{i}^{n-1}\right]$,
and then
\begin{align*}
H &=\left[\begin{array}{cccc}
\sum_{i\in\I}\bar{\omega}_{i} & \sum_{i\in\I}\bar{\omega}_{i}\lambda_{i} & \cdots & \sum_{i\in\I}\bar{\omega}_{i}\lambda_{i}^{n-1}\\
\sum_{i\in\I}\bar{\omega}_{i}\lambda_{i} & \sum_{i\in\I}\bar{\omega}_{i}\lambda_{i}^{2} & \cdots & \sum_{i\in\I}\bar{\omega}_{i}\lambda_{i}^{n}\\
\vdots & \vdots & \ddots & \vdots\\
\sum_{i\in\I}\bar{\omega}_{i}\lambda_{i}^{n-1} & \sum_{i\in\I}\bar{\omega}_{i}\lambda_{i}^{n} & \cdots & \sum_{i\in\I}\bar{\omega}_{i}\lambda_{i}^{2n-2}
\end{array}\right] \\
    &=\sum_{i\in\I}\bar{\omega}_{i}\mathbf{v}_{i}\mathbf{v}_{i}^{\intercal}.
\end{align*}
Post-multiplying $H$ by an arbitrary vector $\mathbf{z}\in \mathbb{R}^n$, we have that $H\mathbf{z}=\sum_{i\in\I}\kappa_i\mathbf{v}_{i}$ where $\kappa_i \coloneqq \bar{\omega}_{i}\mathbf{v}_{i}^{\intercal}\mathbf{z}$. Hence, the column space of $H$ is equal to the span of $\{\mathbf{v}_i\}_{i\in\I}$. Since the $\lambda_{i}$ for which $i\in\I$ are distinct, the $\mathbf{v}_{i}$ are linearly independent. Therefore, the rank of $H$ is equal to $|\I|=|\S_{\mu_{\tL}}|$. \qed

\end{pf}

With this Lemma in hand, we present our main result on recovering the (observable) eigenvalues of the Laplacian.

\begin{thm}\label{thm:id}
Given the sequence of observations $\left(y\left[k\right]\right)_{k=0}^{2n-1}$ from the system in \eqref{eq:dt_integrator}, define the following Hankel matrix
\begin{align}\label{eq:hankel_Y}
Y \coloneqq \left[\begin{array}{cccc}
y[0] & y[1] & \cdots & y[n-1]\\
y[1] & y[2] & \cdots & y[n]\\
\vdots & \vdots & \ddots & \vdots\\
y[n-1] & y[n] & \cdots & y[2n-2]
\end{array}\right],
\end{align}
and denote its rank by $r$. Then, the eigenvalues of $\tL$ which are
in the support of $\mu_{\tL}$ are roots of the polynomial
\[
p_{\tL}\left(x\right)=x^{r}+\alpha_{r-1}x^{r-1}+\cdots+\alpha_{1}x+\alpha_{0},
\]
where the coefficients $\alpha_{0},\ldots,\alpha_{r-1}$ are given
by
\[
\left[\!\begin{array}{c}
\alpha_{0}\\
\alpha_{1}\\
\vdots\\
\alpha_{r-1}
\end{array}\!\right]\!=\!-\!\left[\begin{array}{cccc}
y[0] & y[1] & \cdots & y[r-1]\\
y[1] & y[2] & \cdots & y[r]\\
\vdots & \vdots & \ddots & \vdots\\
y[r-1] & y[r] & \cdots & y[2r-2]
\end{array}\right]^{-1}\!\left[\begin{array}{c}
y[r]\\
y[r+1]\\
\vdots\\
y[2r-1]
\end{array}\right]\!.
\]
\end{thm}

\begin{pf}
By Lemma~\ref{lem:conds}, we know that at most we may recover all eigenvalues $\lambda_i\in\S_{\mu_{\tL}}$. As in Lemma~\ref{lem:rank}, let $\I = \left\{ i\in\{1,\ldots,n\} : \lambda_i \in S_{\mu_{\tL}}\right\}$. By Lemma \ref{lem:rank}, we know that $\text{rk}\left(H\right)=\left|\I\right|$, which we denote by $r$.
For simplicity of exposition, we re-index the $\lambda_{i}$ so that $\I=\{1,\ldots,r\}$.
Define the following polynomial:
\[
p_{\tL}\left(x\right) \coloneqq \prod_{i\in\I}(x-\lambda_i) =  x^{r}+\alpha_{r-1}x^{r-1}+\cdots+\alpha_{1}x+\alpha_{0}.
\]
Notice that, since the eigenvalues are unknown, the coefficients of the polynomial are also unknown. In what follows, we propose an efficient technique to find these coefficients. Since the $\lambda_{i}\in\S_{\mu_{\tL}}$ are roots of $p_{\tL}$, we have the following system of equations:
\begin{align*}
\lambda_{i}^{r}+\alpha_{r-1}\lambda_{i}^{r-1}+\cdots+\alpha_{1}\lambda_{i}+\alpha_{0} & =0, \quad \forall i \in \I.
\end{align*}
Multiplying each equation by the corresponding $\omega_{i}\lambda_i^s$, we obtain
\begin{align*}
\omega_{i}\lambda_{i}^s\left(\lambda_{i}^{r}+\alpha_{r-1}\lambda_{i}^{r-1}+\cdots+\alpha_{1}\lambda_{i}+\alpha_{0}\right) & =0, \quad \forall i\in\I,
\end{align*}
Summing all the equations above over $\I$, we obtain
\begin{align*}
&\sum_{i\in\I} \omega_i\lambda_i^{r+s} + \alpha_{r-1}\sum_{i\in\I} \omega_i\lambda_i^{r-1+s} + \cdots + \alpha_0\sum_{i\in\I} \omega_i\lambda_i^{s} \\
&\quad=m_{r+s}+\alpha_{r-1}m_{r-1+s}+\cdots+\alpha_{1}m_{1+s}+\alpha_{0}m_{s} \\
&\quad=y[r\!+\!s]+\alpha_{r-1}y[r\!-\!1\!+\!s]+\cdots+\alpha_{1}y[1\!+\!s]+\alpha_{0}y[s] \\
&\quad=0,
\end{align*}
where the first equality comes from the definition of $m_k$ in \eqref{eq:moments} and the second comes from the fact that $m_k=y[k]$ for all $k$. Considering the equations obtained for $s \in\{0,1,\ldots,r-1\}$, we obtain a set of linear equations that can be written in matrix form as follows:
\[
\left[\begin{array}{cccc}
y[0] & y[1] & \cdots & y[r-1]\\
y[1] & y[2] & \cdots & y[r]\\
\vdots & \vdots & \ddots & \vdots\\
y[r-1] & y[r] & \cdots & y[2r-2]
\end{array}\right]\left[\begin{array}{c}
\alpha_{0}\\
\alpha_{1}\\
\vdots\\
\alpha_{r-1}
\end{array}\right]=-\left[\begin{array}{c}
y[r]\\
y[r+1]\\
\vdots\\
y[2r-1]
\end{array}\right].
\]
Since $\text{rk}\left(H_{r}\right)=r$ and $m_k = y[k]$, we can find a solution by a simple matrix inversion.
Using the coefficients $\{\alpha_0,\ldots,\alpha_{r-1}\}$, we can compute the roots of $p_{\tL}$ to recover the eigenvalues of $\tL$ that are in the support of $\mu_{\tL}$, i.e., those eigenvalues $\lambda_i$ corresponding to the observable eigenmodes of the Laplacian dynamics. \qed
\end{pf}

While Theorem~\ref{thm:id} makes use of $2n$ observations $\left(y\left[k\right]\right)_{k=0}^{2n-1}$, in practice, fewer observations may be required. Since at most $r = |\S_{\mu_{\tL}}|$ eigenvalues can be recovered, we can build a $k\times k$ Hankel matrix of observations using the first $2k\leq 2r$ observations from the system. Then, we should stop taking observations whenever the rank of this Hankel matrix ceases to grow (i.e., when $k=r$), or when $2n$ observations are obtained, whichever occurs first. In other words, at most $2n$ observations are required to recover the eigenvalues of $\tL$ which correspond to the observable modes of the dynamics, but in practice fewer may be used.

\if{false}
Remark: Under what conditions is our polynomial the characteristic
polynomial of $G$...

Special case $\mathbf{c}=\mathbf{e}_{i}$... TBW...

Special case $\mathbf{c}=\sum_{i\in S\subseteq V}\mathbf{e}_{i}$...
TBW...

Special case $G=G^{\intercal}$, for example, the adjacency or the
Laplacian matrices of an undirected (possibly weighted) graph. In
this case, $\mathbf{u}_{i}=\mathbf{w}_{i}$; hence,
\[
\omega_{j}=\mathbf{c}^{\intercal}\mathbf{u}_{j}\mathbf{u}_{j}^{\intercal}\mathbf{c}=\left(\mathbf{c}^{\intercal}\mathbf{u}_{j}\right)^{2}\geq0.
\]
Furthermore, if $\mathbf{c}=\mathbf{e}_{i}$, then
\[
\omega_{j}=\left(\mathbf{e}_{i}^{\intercal}\mathbf{u}_{j}\right)^{2}=u_{j,i}^{2},
\]
where $u_{j,i}$ is the $i$-th component of the $j$-th eigenvector.
Similarly, for $\mathbf{c}=\sum_{i\in S\subseteq V}\mathbf{e}_{i}$,
we have,
\[
\omega_{j}=\left(\sum_{i\in S\subseteq V}\mathbf{e}_{i}^{\intercal}\mathbf{u}_{j}\right)^{2}=\left(\sum_{i\in S\subseteq V}u_{j,i}\right)^{2}\geq0.
\]

\fi 

\subsection{Network of DT Identical Agents}\label{subsec:dt_multi}

In many applications, the network of interest will not only contain single integrators, but instead will consist of agents with more general dynamics. With this in mind, consider a network of $n$ agents where each agent follows the dynamics
$\mathbf{x}_{i}\left[k+1\right]=A\mathbf{x}\left[k\right]+\mathbf{u}\left[k\right]$, where\textbf{ $\mathbf{x}_{i}$} is a $d$-dimensional vector of states
and $\mathbf{u}\left[k\right]$ is a linear combination of the states of
the neighboring agents of $i$, i.e.,
\begin{align}
\begin{split}\label{eq:dt_multi}
\mathbf{x}_{i}\left[k+1\right] &=A\mathbf{x}\left[k\right]+\sum_{j=1}^{n}l_{ij}\mathbf{x}_{j}\left[k\right],\;\mathbf{x}_{i}\left[0\right]=x_{0i}\boldb,\\
y\left[k\right] & =\sum_{i=1}^{n}c_{i}\boldg^{\intercal}\mathbf{x}_{i}\left[k\right],
\end{split}
\end{align}
where $l_{ij} = [\tL]_{ij}$, $c_i = [\mathbf{c}]_i$, $x_{0i} = [\x_0]_i$, and $\boldb$ and $\boldg$ are vectors of the individual states. In other words, all agents start with the initial condition $\boldb$ weighted by $x_{0i}$, and all individual observations are $\boldg^\intercal\x_i[k]$  weighted by $c_i$. Stacking the vectors of states in a large vector $\mathbf{x}=\left(\mathbf{x}_{1}^{\intercal},\ldots,\mathbf{x}_{n}^{\intercal}\right)^{\intercal}$,
the dynamics can be written as
\begin{align*}
\mathbf{x}\left[k+1\right] & =\left(I_{n}\otimes A+\tL\otimes I_{d}\right)\mathbf{x}\left[k\right],\;\mathbf{x}\left[0\right]=\mathbf{x}_0\otimes\boldb,\\
y\left[k\right] & =\left(\mathbf{c}\otimes\boldg\right)^{\intercal}\mathbf{x}\left[k\right].
\end{align*}
We assume the state matrix $A$ of each agent is known, but the
(normalized) Laplacian matrix $\tL$ is unknown. We aim towards reconstructing
the (observable) eigenvalue spectrum of $\tL$ from a finite sequence of outputs. A technical difficulty in this case is that we no longer have $m_k = y[k]$; however, as we will see, we may still recover the moments using the finite sequence of observations. This result is summarized in the following theorem.

\begin{thm}\label{thm:id_multi}
Given the sequence of observations $\left(y\left[k\right]\right)_{k=0}^{2n-1}$ from the system in \eqref{eq:dt_multi}, consider the Hankel matrix $Y$ defined in \eqref{eq:hankel_Y} and denote its rank by $r$. The moments of $\mu_{\tL}$ satisfy the following equality:
\begin{align*}
\small{
\!\begin{bmatrix}
m_{0}\\
m_{1}\\
\vdots\\
\!m_{2r-1}\!
\end{bmatrix} \!=\! \begin{bmatrix}
b_{0,0}\nu_{0} & 0 & \cdots & 0\\
b_{1,0}\nu_{1} & b_{1,1}\nu_{0} & \cdots & 0\\
\vdots & \vdots & \ddots & \vdots\\
\!b_{2r-1,0}\nu_{2r-1} \!&\! b_{2r-1,1}\nu_{2r-2} \!&\! \cdots \!&\! b_{2r-1,2r-1}\nu_{0}\!
\end{bmatrix}^{\!-1\!}\!
\begin{bmatrix}
y_{0}\\
y_{1}\\
\vdots\\
\!y_{2r-1}\!
\end{bmatrix}\!
}
\end{align*}
where $\nu_{k-s} \coloneqq \boldg^{\intercal}A^{k-s}\boldb$, $b_{k,s} \coloneqq \binom{k}{s}$, and the matrix is invertible when $\boldg^\intercal\boldb \neq 0$. Then, the eigenvalues of $\tL$ contained in the support of $\mu_{\tL}$ are roots of the polynomial
\[
p_{\tL}\left(x\right)=x^{r}+\alpha_{r-1}x^{r-1}+\cdots+\alpha_{1}x+\alpha_{0},
\]
where the coefficients $\alpha_{0},\ldots,\alpha_{r-1}$ satisfy
\[
\left[\!\begin{array}{c}
\alpha_{0}\\
\alpha_{1}\\
\vdots\\
\alpha_{r-1}
\end{array}\!\right]\!=\!-\!\left[\begin{array}{cccc}
m_0 & m_1 & \cdots & m_{r-1}\\
m_1 & m_2 & \cdots & m_r\\
\vdots & \vdots & \ddots & \vdots\\
m_{r-1} & m_r & \cdots & m_{2r-2}
\end{array}\right]^{-1}\!\left[\begin{array}{c}
m_r\\
m_{r+1}\\
\vdots\\
m_{2r-1}
\end{array}\right]\!.
\]
\end{thm}

\begin{pf}
Considering the diagonalization $\tL=U\Lambda U^{-1}$, we have
\begin{align*}
&\left(I_{n}\otimes A+\tL\otimes I_{d}\right)^{k} \\
 &~= \left[\left(U\otimes I_{d}\right)\left(I_{n}\otimes A+\Lambda\otimes I_{d}\right)\left(U^{-1}\otimes I_{d}\right)\right]^k \\
 &~=\left(U\otimes I_{d}\right)\left(I_{n}\otimes A+\Lambda\otimes I_{d}\right)^{k}\left(U^{-1}\otimes I_{d}\right)\\
 &~= \left(U\otimes I_{d}\right)\left[\sum_{s=0}^{k}\binom{k}{s}\left(I_{n}\otimes A^{k-s}\right)\left(\Lambda^{s}\otimes I_{d}\right)\right]\left(U^{-1}\otimes I_{d}\right).
\end{align*}
Thus,
\begin{align*}
&y\left[k\right] = \left(\mathbf{c}\otimes\boldg\right)^{\intercal}\left(I_{n}\otimes A+\tL\otimes I_{d}\right)^{k}\left(\mathbf{x}_0\otimes\boldb\right) \\
&=\sum_{s=0}^{k}\!\binom{k}{s}\! \left(\mathbf{c}^{\intercal}U\otimes\boldg^{\intercal}\right) \! \left(I_n\otimes A^{k-s}\right) \! \left(\Lambda^{s}\otimes I_{d}\right) \! \left(U^{-1}\mathbf{x}_0\otimes\boldb\right)\\
 & =\sum_{s=0}^{k}\binom{k}{s}\left(\mathbf{c}^{\intercal}U\Lambda^{s}U^{-1}\mathbf{x}_0\right)\left(\boldg^{\intercal}A^{k-s}\boldb\right).
\end{align*}
Hence, we obtain
\begin{align}\label{eq:obs_multi}
y\left[k\right]=\sum_{s=0}^{k}\binom{k}{s}\nu_{k-s}\sum_{i=1}^{n}\omega_{i}\lambda_{i}^{s}=\sum_{s=0}^{k}\binom{k}{s}\nu_{k-s}m_{s}.
\end{align}
From the sequence $\left(y\left[k\right]\right)_{k=0}^{2n-1}$, we obtain a lower triangular system of linear equations that can be solved to find the sequence of moments $\left(m_{k}\right)_{k=0}^{2n-1}$. Specifically, if we collect $2r$ observations, with $b_{k,s} = \binom{k}{s}$, we have that \eqref{eq:obs_multi} for $k=0,\ldots,2r-1$ results in
\begin{align*}
\small{
\!\begin{bmatrix}
y_{0}\\
y_{1}\\
\vdots\\
\!y_{2r-1}\!
\end{bmatrix} \!=\! \begin{bmatrix}
b_{0,0}\nu_{0} & 0 & \cdots & 0\\
b_{1,0}\nu_{1} & b_{1,1}\nu_{0} & \cdots & 0\\
\vdots & \vdots & \ddots & \vdots\\
\!b_{2r-1,0}\nu_{2r-1} & b_{2r-1,1}\nu_{2r-2} \!&\! \cdots \!&\! b_{2r-1,2r-1}\nu_{0}\!
\end{bmatrix}\begin{bmatrix}
m_{0}\\
m_{1}\\
\vdots\\
\!m_{2r-1}\!
\end{bmatrix}
}
\end{align*}
As long as $\nu_0 = \boldg^\intercal\boldb \neq 0$, the above matrix is full-rank. We may then recover the moments by a simple inversion, and apply Theorem~\ref{thm:id} to find the eigenvalues of $\tL$. \qed
\end{pf}

\section{Continuous-Time Dynamics}\label{sec:ct}

In the case of continuous-time (CT) dynamics, there are some subtle but important differences to the case of discrete-time. Fortunately, our main results are still applicable in this domain, as we will describe in the following subsections.

\begin{figure*}[!ht]
\centering
    \begin{subfigure}{.32\textwidth}
    \centering
    \includegraphics[width=\linewidth]{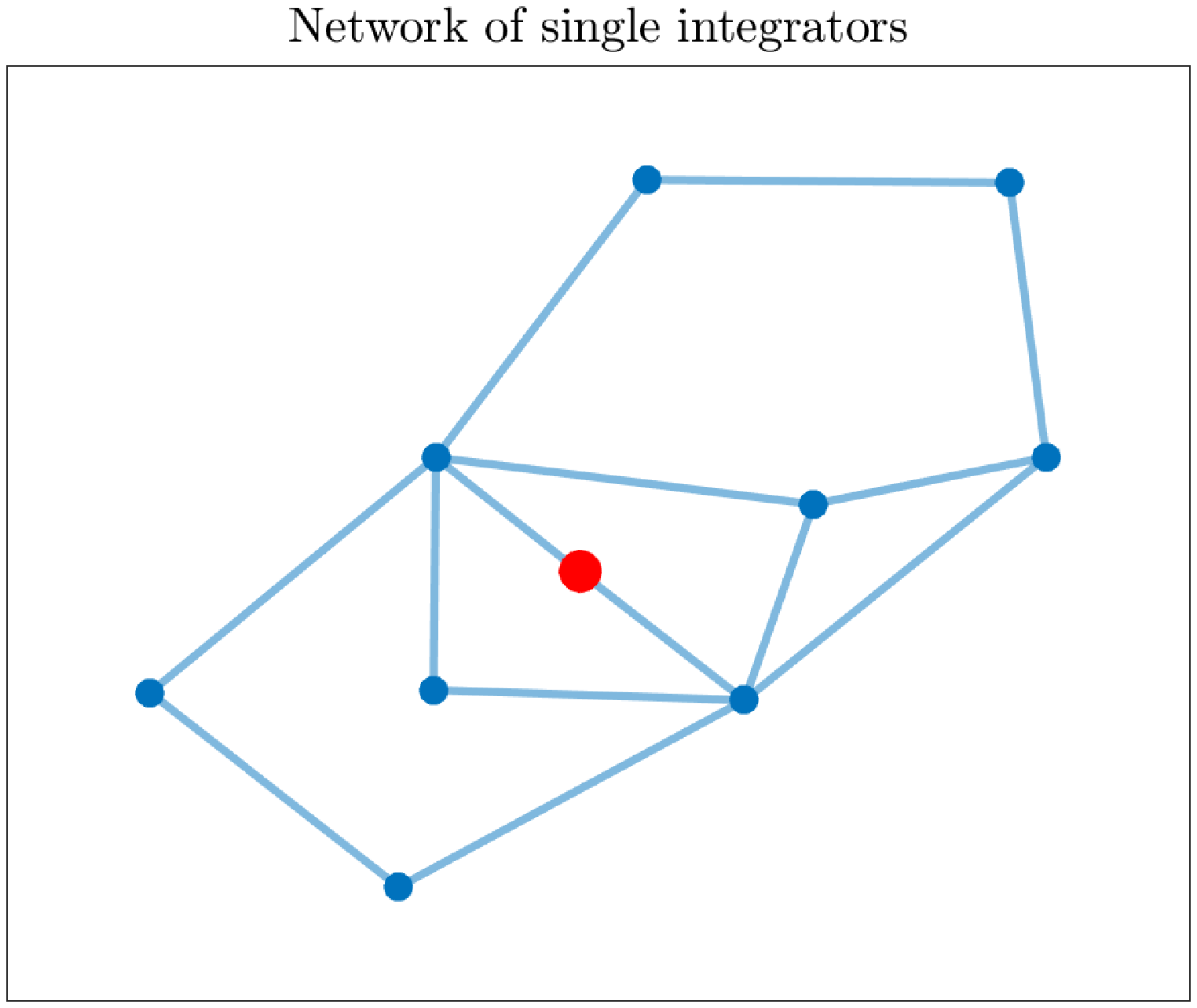}
    \caption{Network topology, with single output agent highlighted.}
    \label{fig:dt_ntwk}
    \end{subfigure}
    \begin{subfigure}{.3\textwidth}
    \centering
    \includegraphics[width=\linewidth]{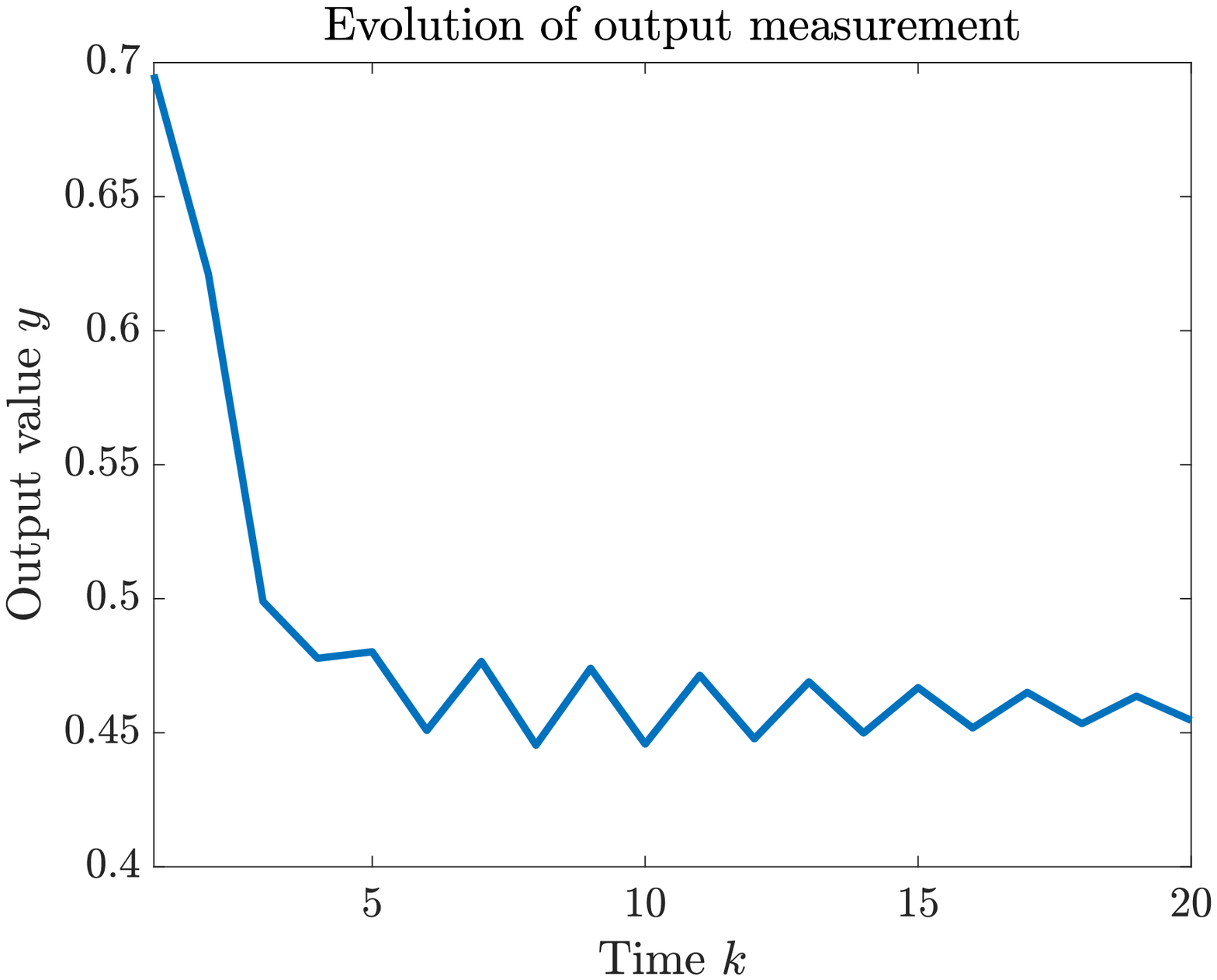}
    \caption{Output $y[k] = e_i^\intercal\tL x_0$, where we observe only agent $i$.}
    \label{fig:dt_y}
    \end{subfigure}
    \begin{subfigure}{.3\textwidth}
    \centering
    \includegraphics[width=\linewidth]{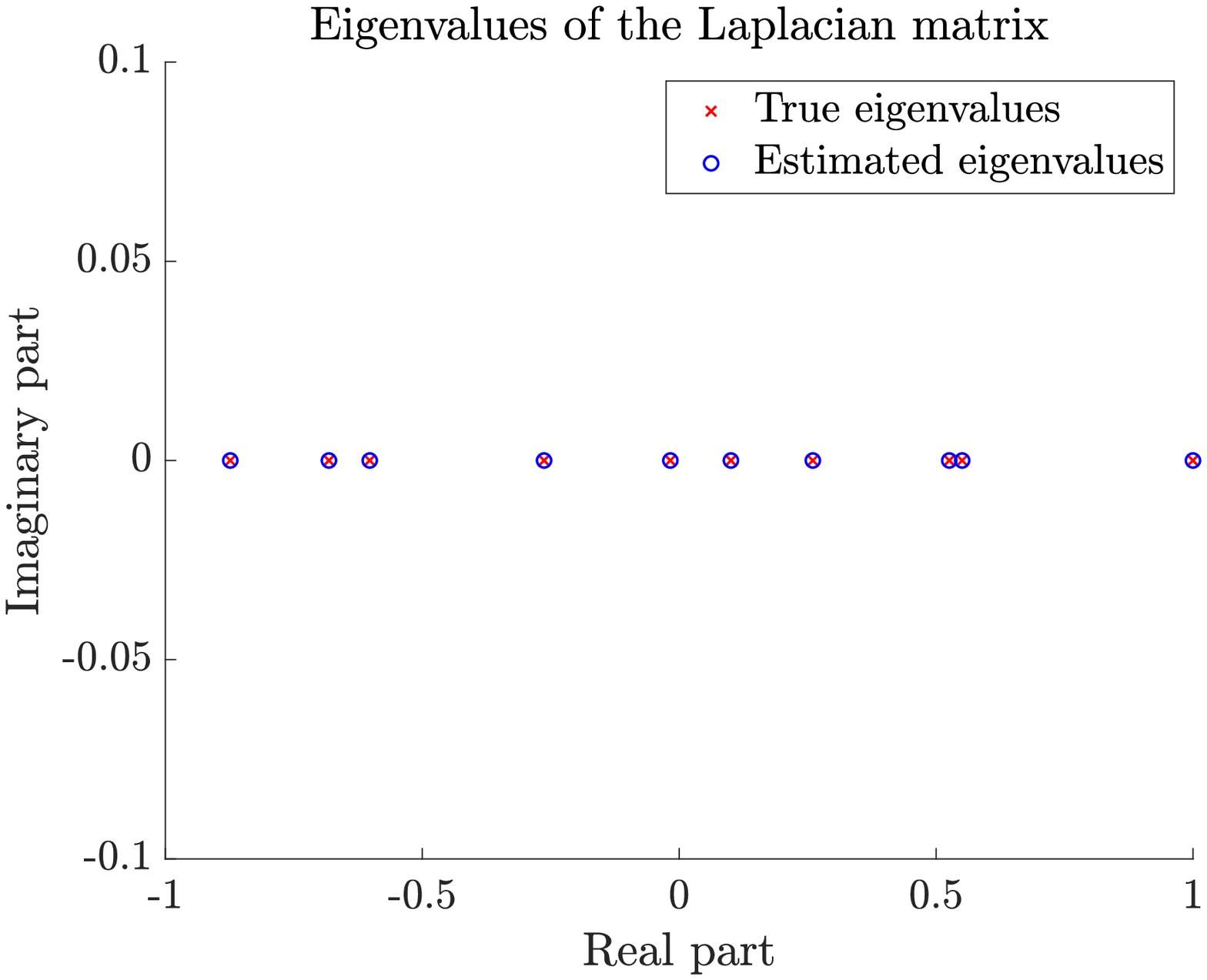}
    \caption{Comparison of true and estimated eigenvalues.}
    \label{fig:dt_eigs}
    \end{subfigure}
    \caption{$10$-agent preferential attachment network in discrete-time, generated according to \cite{Barabasi99}. The initial condition is randomly generated as $\x_0 \sim \text{Uniform}[0,1]^n$. There are $10$ unique eigenvalues of $\tL$ in this case, which are all recovered via our estimation approach.}
    \label{fig:dt_sims}
\end{figure*}
\begin{figure*}[!ht]
\centering
    \begin{subfigure}{.32\textwidth}
    \centering
    \includegraphics[width=\linewidth]{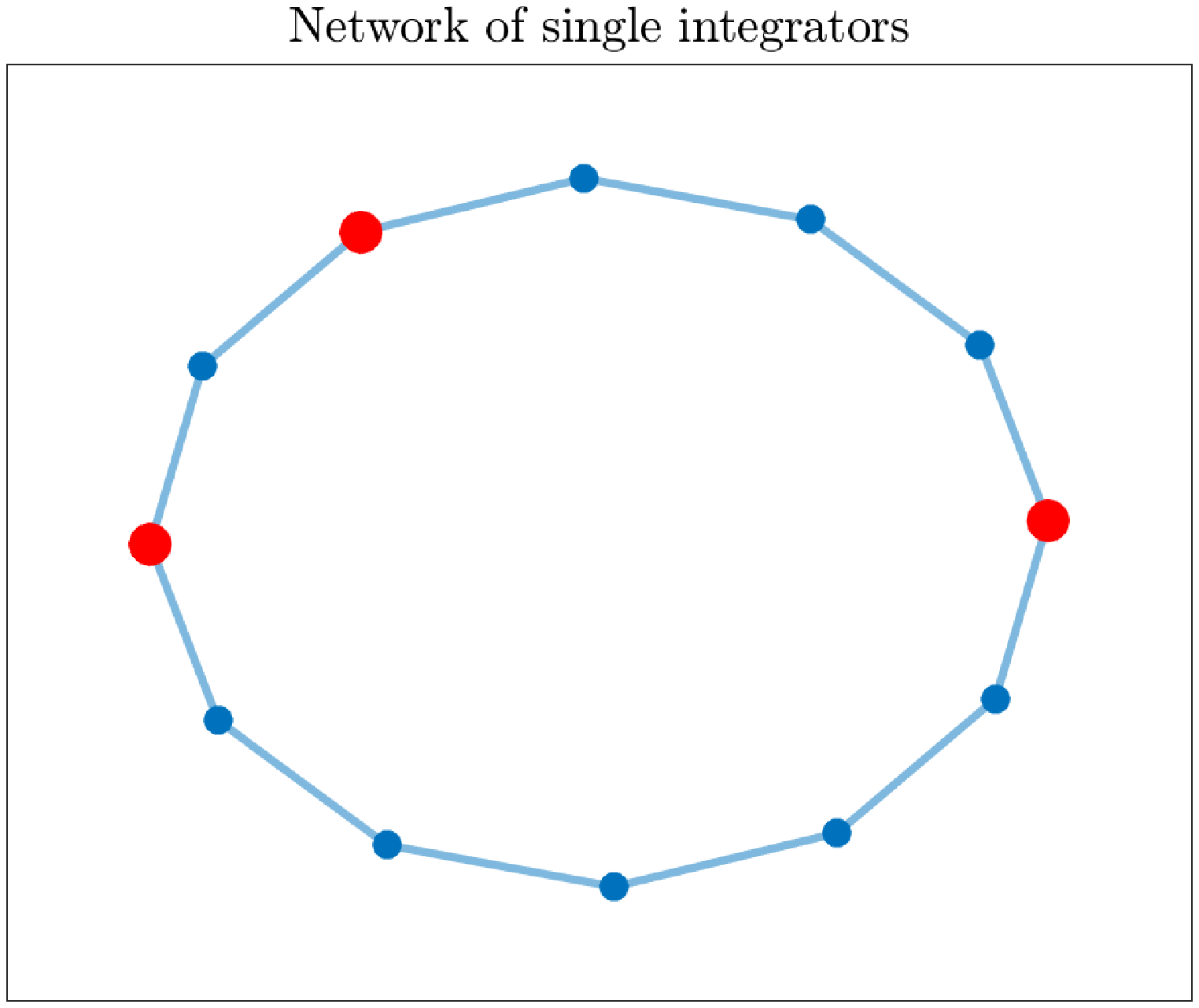}
    \caption{Network topology, with output agents highlighted.}
    \label{fig:ct_ntwk}
    \end{subfigure}
    \begin{subfigure}{.3\textwidth}
    \centering
    \includegraphics[width=\linewidth]{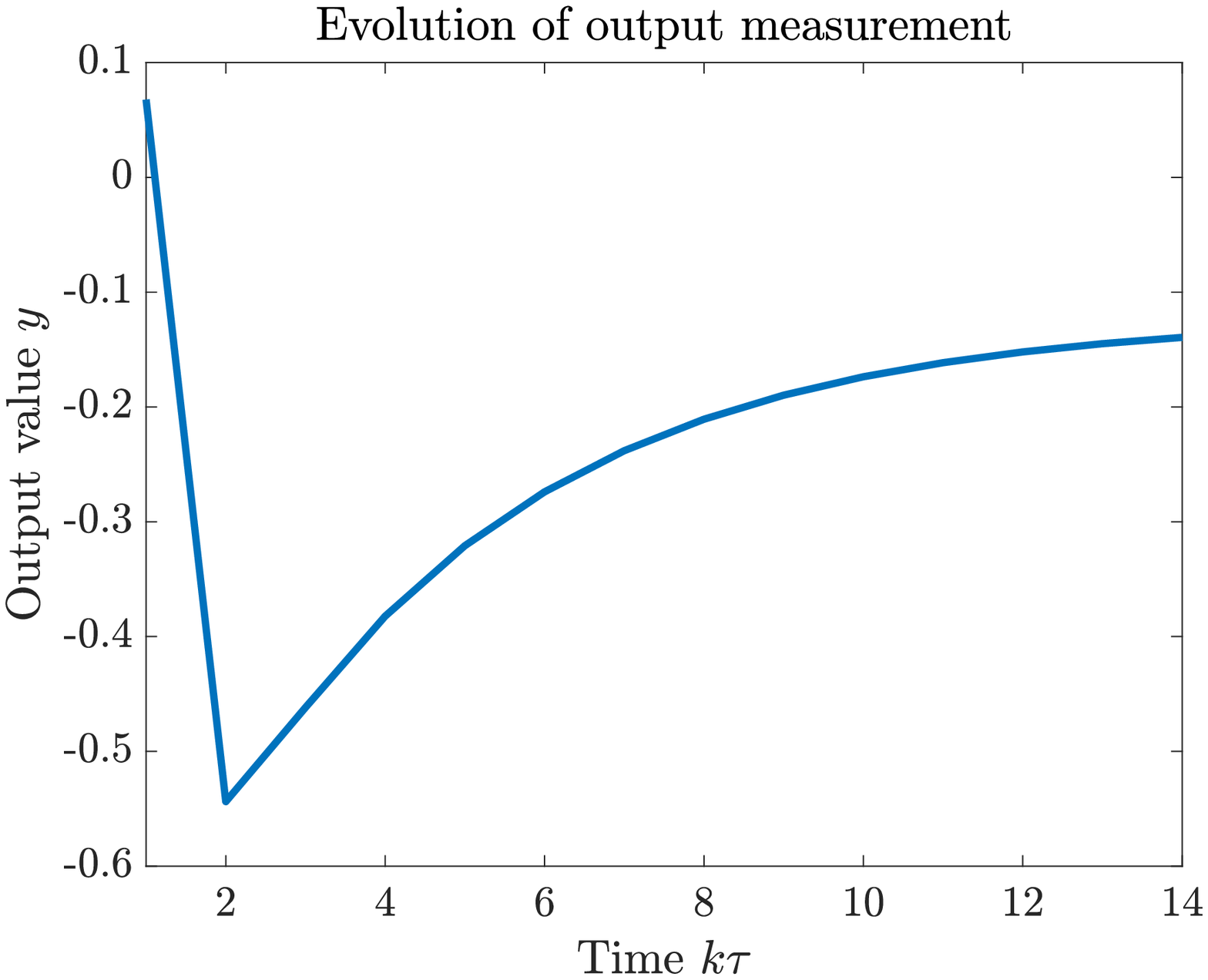}
    \caption{Output $y[k] = c^\intercal e^{-Lk\tau} x_0$, where we observe three agents with equal weight.}
    \label{fig:ct_y}
    \end{subfigure}
    \begin{subfigure}{.3\textwidth}
    \centering
    \includegraphics[width=\linewidth]{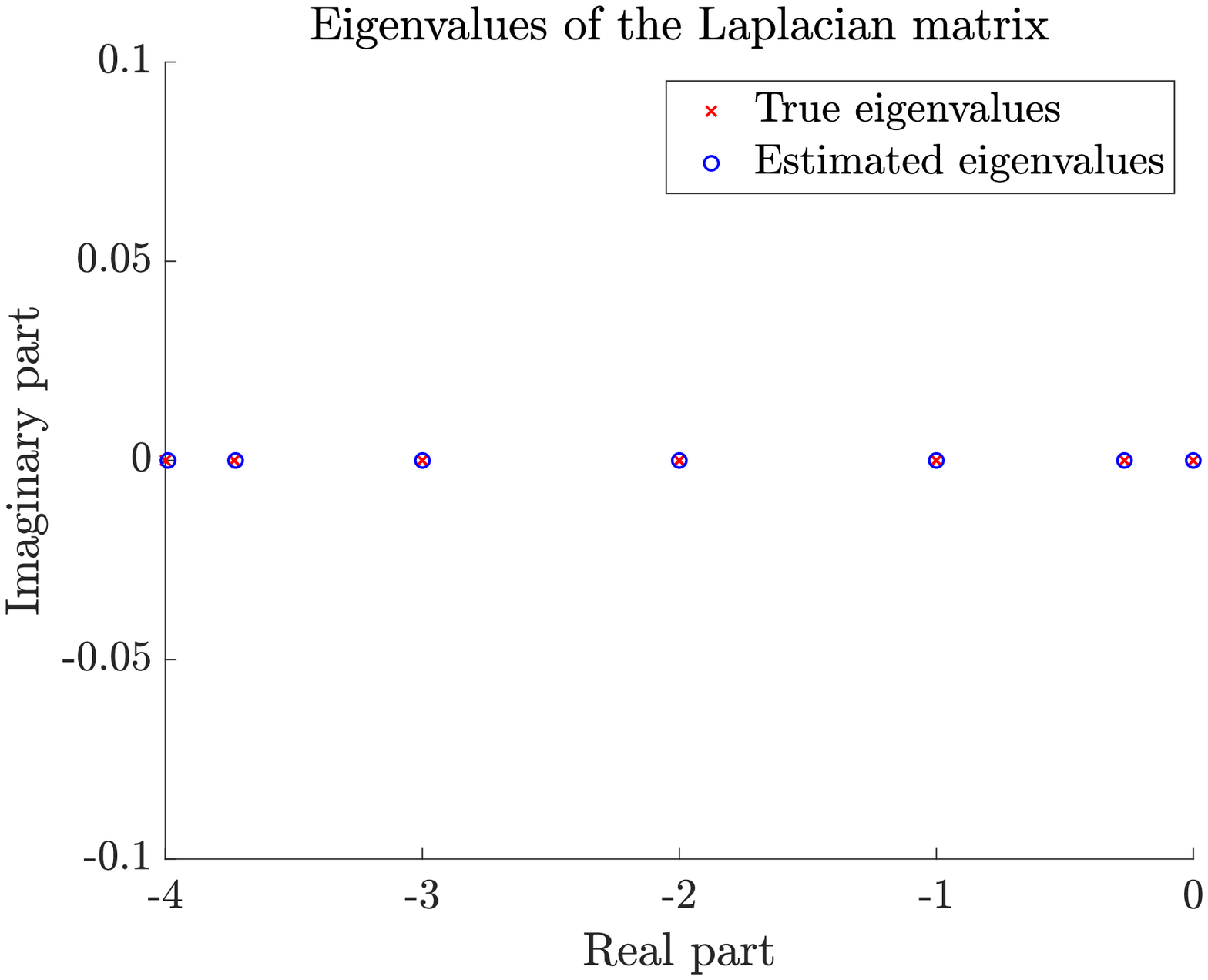}
    \caption{Comparison of true and estimated eigenvalues (repeated values are overlaid).}
    \label{fig:ct_eigs}
    \end{subfigure}
    \caption{$12$-agent single integrator ring network in continuous-time, with sampling rate $\tau = 1$ and random initial condition $\x_0 \sim \text{Uniform}[0,1]^n$. Here there are $7$ unique eigenvalues of $L$, all of which are recovered via our estimation approach.}
    \label{fig:ct_sims}
\end{figure*}

\subsection{Network of Single Integrators}\label{subsec:ct_single}

Consider a network of coupled CT single integrators:
\begin{align*}
\dot{\mathbf{x}}(t) & = -L\mathbf{x}(t),\;\mathbf{x}\left(0\right)=\mathbf{x}_0,\\
y(t) & =\mathbf{c}^{\intercal}\mathbf{x}(t),
\end{align*}
where $L \coloneqq D - G$ is the combinatorial Laplacian, which is symmetric and, hence, diagonalizable with real eigenvalues as $L = U\Lambda U^\intercal$. We thus have
\[
y\left(t\right)=\mathbf{c}^{\intercal}e^{-Lt}\mathbf{x}_0.
\]
In practice, we take discrete samples of the output with an arbitrary period $\tau>0$; diagonalizing $L=U\Lambda W$, we have
\begin{align*}
y_{k} &\coloneqq y\left(k\tau\right)
= \mathbf{c}^{\intercal}U e^{-\Lambda k\tau}W\mathbf{x}_0 \\
&= \sum_{i=1}^{n}\omega_{i}\left(e^{-\lambda_{i}\tau}\right)^{k},
\end{align*}
since $\Lambda$ is diagonal and $W = U^\intercal = U^{-1}$, with $\omega_{i}$ as defined in \eqref{eq:weights}.

Similarly to the discrete-time case, we define the following signed Borel measure on $\mathbb{R}:$
\[
\mu_{L}\left(z\right) \coloneqq \sum_{i=1}^{n}\omega_{i}\delta\left(z-e^{-\lambda_{i}\tau}\right),
\]
with $k$-th moment given by
\begin{align*}
m_{k} &\coloneqq \int_{\mathbb{R}} z^{k}\sum_{i=1}^{n}\omega_{i}\delta\left(z-e^{-\lambda_{i}\tau}\right)dz=\sum_{i=1}^{n}\omega_{i}\left(e^{-\lambda_{i}\tau}\right)^{k} = y_k.
\end{align*}
In contrast with the discrete-time case, the support of this measure is $S_{\mu_L} \coloneqq \left\{ e^{-\lambda_{i}\tau}\colon\sum_{j\colon\lambda_{j}=\lambda_{i}}\omega_{j}\neq0\right\} $. However, we may still apply Theorem \ref{thm:id} to recover the support of $\mu_{L}$, i.e., the quantities $e^{-\lambda_i\tau}$, from the finite sequence $\left(y_{k}\right)_{k=0}^{2n-1}$. The eigenvalues of the Laplacian matrix corresponding to the observable modes of the system may then be recovered by taking a logarithm and dividing by~$-\tau$.

\subsection{Network of Identical CT Agents}\label{subsec:ct_multi}

Similarly to the more general setting in Section~\ref{subsec:dt_multi}, we consider the dynamics of a network of CT agents, which can be described (in a compact form) as follows 
\begin{align*}
\dot{\mathbf{x}}(t) & =\left(I_{n}\otimes A - L\otimes I_{d}\right)\mathbf{x}(t),\;\mathbf{x}\left(0\right)=\mathbf{x}_0\otimes\boldb,\\
y(t) & =\left(\mathbf{c}\otimes\boldg\right)^{\intercal}\mathbf{x}(t).
\end{align*}
Hence, considering a sampling period $\tau>0$, we have that
\begin{align*}
y\left(k\tau\right) & =\left(\mathbf{c}\otimes\boldg\right)^{\intercal}e^{\left(I_{n}\otimes A - L\otimes I_{d}\right)k\tau}\left(\mathbf{x}_0\otimes\boldb\right)\\
 & =\left(\mathbf{c}\otimes\boldg\right)^{\intercal}\left(e^{-Lk\tau}\otimes e^{Ak\tau}\right)\left(\mathbf{x}_0\otimes\boldb\right)\\
 & =\left(\mathbf{c}^{\intercal}U e^{-\Lambda k\tau}W\mathbf{x}_0\right)\left(\boldg^{\intercal}e^{Ak\tau}\boldb\right)\\
 & =\nu_{k}\sum_{i=1}^{n}\omega_{i}\left(e^{-\lambda_{i}\tau}\right)^{k},
\end{align*}
where $\nu_{k}\coloneqq\boldg^{\intercal}e^{Ak\tau}\boldb$ and $\omega_{i}$ is defined in (\ref{eq:weights}). Applying Theorem~\ref{thm:id_multi} followed by a logarithmic transformation, we again obtain the eigenvalues of the Laplacian matrix.

\section{Simulations}\label{sec:sims}

In this section we illustrate our results on a variety of simulated networks. In each case, some underlying network structure is created which is unknown to us. We then simulate the evolution of the system with random initial condition $\x_0$ and an observability vector $\c$ which is unknown to the algorithm, and then compare the estimated eigenvalues to the true spectrum of the Laplacian matrix.

Figure~\ref{fig:dt_sims} shows the result of using Theorem~\ref{thm:id} on the randomly generated $10$-agent preferential attachment network shown in Figure~\ref{fig:dt_ntwk} (see \cite{Barabasi99}), where we model each agent using a single integrator dynamics in discrete-time, as in~\eqref{eq:dt_integrator}. In Figure~\ref{fig:dt_y}, we show the evolution of our output signal measured from a single agent highlighted in the figure. Figure~\ref{fig:dt_eigs} compares both the true and estimated eigenvalues. In this case there are $10$ unique eigenvalues of $\tL$ and all of these are perfectly recovered using a sequence of 20 measurements retrieved from a single agent.

In Figure~\ref{fig:ct_sims} we see the results of using Theorem~\ref{thm:id} followed by a logarithmic transformation on the $12$-agent ring network shown in Figure~\ref{fig:ct_ntwk} using the continuous-time single integrator dynamics from Section~\ref{sec:ct}. In this case, we observe the sum of the measurements obtained from the three highlighted agents; the overall measurement is shown in Figure~\ref{fig:ct_y}. This network has only $7$ unique eigenvalues, all of which are estimated (without multiplicities) by our algorithm after recovering the support of $\S_{\mu_L}$ and performing a logarithmic transformation, as shown in Figure~\ref{fig:ct_eigs}.

\section{Conclusion}\label{sec:conc}

In this paper, we have proposed an efficient methodology for recovering the observable eigenvalues of the Laplacian matrix of a network of interacting dynamical agents using a sparse set of output measurements. Unlike other methods, we require only a single finite sequence of measurements from the multiagent network of length, at most, $2n$. Moreover, we need no prior knowledge of the network topology, initial condition, or which agents are contributing to the measurements. We develop our technique for systems in both discrete- and continuous-time. We consider the case of agents modeled by single integrators, as well as more complex dynamics. Our simulation results show that we are able to recover the spectrum of the Laplacian in all cases with high accuracy.

\bibliography{MH-bib}

\end{document}